\documentclass[11pt]{article}
\usepackage{amsmath,amsfonts}
\usepackage{txfonts,wasysym,lscape,amssymb,mathrsfs,extarrows}
\usepackage[all,pdf]{xy}
\RequirePackage{color,graphicx}

\newtheorem{theorem}{Theorem}[section]
\newtheorem{proposition}{Proposition}[section]
\newtheorem{lemma}{Lemma}[section]

\newtheorem{definition}{Definition}[section]

\title{Stability for Constrained Minimax Optimization}

\author{Yu-Hong Dai\footnote{LSEC, ICMSEC, AMSS, Chinese Academy of Sciences, Beijing 100190, China.  {\sl Email}: dyh@lsec.cc.ac.cn.
This author was supported by the Natural Science Foundation of China (Nos. 11991020, 11631013, 11971372 and 11991021) and the Strategic Priority Research Program of Chinese Academy of Sciences (No. XDA27000000).}
 \footnote{School of Mathematical Sciences, University of Chinese Academy of Sciences, Beijing 100049, China.}\quad and \quad Liwei Zhang \footnote{Corresponding author. School of Mathematical Sciences, Dalian University of Technology, Dalian 116024, China. {\sl Email}: lwzhang@dlut.edu.cn. This author was supported by the Natural Science Foundation of China (Nos. 11971089 and 11731013).}}

\date{}
\begin{document}

\maketitle

\begin{abstract}
Minimax optimization problems are an important class of optimization problems arising from both modern machine learning and from traditional research
areas. We focus on  the stability of constrained minimax optimization problems based on the notion of local minimax point by Dai and Zhang (2020). Firstly, we extend the  classical Jacobian uniqueness conditions of nonlinear programming to the constrained minimax problem and prove that this set of properties  is stable with respect to small ${\cal C}^2$ perturbation.  Secondly, we provide a set of conditions, called {\sc Property A}, which does not require the strict complementarity condition for the upper level constraints. Finally, we prove that {\sc Property A} is a sufficient condition for the strong regularity of  the Kurash-Kuhn-Tucker (KKT) system at  the KKT point, and it is also a sufficient condition for
the  local Lipschitzian homeomorphism of the Kojima mapping near the KKT point.

\vskip 6 true pt \noindent \textbf{Key words}: constrained minimax optimization,  Jacobian uniqueness conditions, strong regularity, strong sufficient optimality condition, Kojima mapping, local Lipschitzian homeomorphism.
\vskip 12 true pt \noindent \textbf{AMS subject classification}: 90C30
\end{abstract}
\bigskip\noindent

%=================================================1=============================================
\section{Introduction}
 \setcounter{equation}{0}
\quad \,
Let $m$, $n$, $m_1$, $m_2$, $n_1$ and $n_2$ be positive integers and $f:\Re^n\times \Re^m \rightarrow \Re$, $h:\Re^n\times \Re^m \rightarrow \Re^{m_1}$, $g:\Re^n\times \Re^m \rightarrow \Re^{m_2}$, $H:\Re^n \rightarrow \Re^{n_1}$ and $G:\Re^n \rightarrow \Re^{n_2}$ be given functions. We are
interested in the constrained minimax optimization problem of the form
\begin{equation}\label{cminimax}
\min_{x \in \Phi}\max_{y \in Y(x)}f(x,y),
\end{equation}
where $\Phi\subset \Re^n$ is a feasible set of decision variable $x$ defined by
\begin{equation}\label{Phi}
\Phi=\{x \in \Re^n: H(x)=0,\, G(x)\leq 0\}
\end{equation}
and $Y: \Re^n \rightrightarrows \Re^m$ is a set-valued mapping defined by
\begin{equation}\label{Yx}
Y(x)=\{y \in \Re^m: h(x,y)=0,\, g(x,y)\leq 0\}.
\end{equation}

For unconstrained nonconvex-nonconcave minimax optimization, Jin {\sl et al.} \cite{Jin2019} proposed a
 proper definition of local minimax point. This definition of local minimax point is extended in \cite{DaiZhang2020} for the constrained minimax optimization problem (\ref{cminimax}).
\begin{definition}\label{def:minimaxpoint}
A point $(x^*,y^*) \in \Re^n\times \Re^m$ is said to be a local minimax point of Problem (\ref{cminimax}) if  there exist $\delta_0>0$ and a function $\eta:(0,\delta_0] \rightarrow \Re_+$
 satisfying $\eta(\delta)\rightarrow 0$ as $\delta\rightarrow 0$ such that for any $\delta\in (0, \delta_0]$  and any $(x, y)\in [\textbf{B}_{\delta}(x^*)\cap \Phi] \times [Y(x^*) \cap \textbf{B}_{\delta}(y^*)]$, we have
 \begin{equation}\label{eq:definqs}
 f(x^*,y) \leq f(x^*,y^*)\leq \max_z\left\{f(x,z):z \in Y(x) \cap \textbf{B}_{\eta(\delta)}(y^*)\right\}.
 \end{equation}
\end{definition}
In \cite{DaiZhang2020}, we established the first-order optimality, the second-order necessary and sufficient optimality conditions for Problem (\ref{cminimax}) when the Jacobian uniqueness conditions are satisfied for the lower level problem and the first-order necessary  optimality conditions when the strong second-order sufficient optimality condition and the linear independence constraint qualification are satisfied for the lower level problem.

It is well known that, for nonlinear programming, the Jacobian uniqueness condition can be used to establish the stability of the ${\cal C}^2$-perturbation (see for instance \cite{FiaccoMc90}) and prove that the strong second-order sufficient optimality condition and the linear independence constraint qualification are equivalent to the strong regularity of the Kurash-Kuhn-Tucker (KKT) system (see \cite{Robinson80} and \cite{Jongen90}). The question naturally arises: What are the counterparts of these two stability properties for the constrained minimax optimization problem?  The purpose of this paper is to answer this basic question.

The rest of this paper is organized as follows. In Section 2, we develop a simplified version for second-order optimality conditions for the constrained minimax optimization problem, which is suitable for the study of  stability properties.
In Section 3, we prove that the proposed  Jacobian uniqueness conditions for Problem (\ref{cminimax}) are kept when a ${\cal C}^2$-perturbation of the original problem occurs. In Section 4, we prove that  the proposed
{\sc Property A}, which does not require the strict complementarity for the upper level problem, is a sufficient condition for the strong regularity of  the KKT system at the KKT point. Finally, we draw a conclusion in Section 5.

{\bf Notation}. Scalars and vectors are expressed in lower case
letters and matrices are expressed in upper case letters.
 For a vector $x$, denote $\textbf{B}_{\delta}(x)=\{x':\|x'-x\|\leq \delta\}$. For $a$, $b \in \Re^p$, $a\circ b$ denotes the Hadamard product of $a$ and $b$; namely, $a\circ b=(a_1b_1,\ldots,a_pb_p)^T$. For $a
  \in \Re^p$, $a>0$, denote $\sqrt{a}={\rm Diag}(\sqrt{a_1},\ldots, \sqrt{a_p})$. For a convex set $D \subset \Re^k$, we use $\Pi_D(w)$ to stand for the projection of $w$ onto $D$. For simplicity, for a function $F:\Re^n \times \Re^m \rightarrow \Re$, a mapping $g:\Re^n \times \Re^m \rightarrow \Re^p$,  and a mapping $y:\Re^n\rightarrow \Re^m$, we denote
  $$
  \begin{array}{ll}
  \nabla_xF(x,y(x))=\nabla_x F(x,y)|_{y=y(x)}, &\quad \nabla_yF(x,y(x))=\nabla_y F(x,y)|_{y=y(x)},\\[4pt] \nabla^2_{xx}F(x,y(x))=\nabla^2_{xx} F(x,y)|_{y=y(x)}, &\quad
  \nabla^2_{xy}F(x,y(x))=\nabla^2_{xy} F(x,y)|_{y=y(x)},\\[4pt] {\cal J}_xg(x,y(x))={\cal J}_xg(x,y)|_{y=y(x)}, &\quad {\cal J}_yg(x,y(x))={\cal J}_yg(x,y)|_{y=y(x)}.
  \end{array}
  $$
  Let $G: \Re^n \rightarrow \Re^m$ be a locally Lipschitz continuous mapping over an open set ${\cal O}$. Then $G$ is differentiable almost everywhere in ${\cal O}$. Let ${\cal D}_G$ denote the set of differentiable points of $G$ in ${\cal O}$. For a point $x\in {\cal O}$, the B-subdifferential of $G$ at $x$ is defined by
  $$
  \partial_B G(x)=\left\{V: \exists x^k \in {\cal D}_G, \, x^k \rightarrow x, \, \displaystyle {\cal J}G(x^k) \rightarrow V\right\}
  $$
  and the Clarke subdifferential of $G$ at $x$ is defined by
  $$
  \partial G(x)={\rm conv}\,\partial_B G(x).
  $$
  For differential properties of Lipschitz mappings, see the famous book \cite{Clarke83}.
\section{Simplified Second-order Optimality Conditions}\label{Sec2}
\setcounter{equation}{0}
\quad \, Consider the case when the Jacobian uniqueness conditions hold at some point $(x^*,y^*,\mu^*,\lambda^*)\in \Re^n \times \Re^m \times \Re^q\times \Re^p$, where $(x^*,y^*)\in \Re^n \times \Re^m$ is a point around which  $f,h,g$ are twice continuously differentiable. For a point $x \in \Re^n$ around $x^*$, we use (${\rm P}_x$) to denote the following problem
\begin{equation}\label{eq:Px}
\begin{array}{cl}
\max_{z\in \Re^m} & f(x,z)\\[4pt]
{\rm s.t.\ \ \ \ \ }& h(x,z)=0,\\[4pt]
& g(x,z) \leq 0.
\end{array}
\end{equation}
The Lagrangian of Problem (${\rm P}_x$) is defined by
$$
{\cal L}(x,z,\mu,\lambda)=f(x,z)+\mu^Th(x,z)-\lambda^Tg(x,z).
$$
\begin{definition}\label{def:Jac}
Let $(\mu^*,\lambda^*) \in \Re^{m_1}\times \Re^{m_2}$ be a point. We say that Jacobian uniqueness conditions of Problem
(${\rm P}_{x^*}$) are satisfied at $(y^*,\mu^*,\lambda^*)$ if
\begin{itemize}
\item[{\rm (a)}]The point $(y^*,\mu^*,\lambda^*)$ is a Karush-Kuhn-Tucker point of  Problem
(${\rm P}_{x^*}$); {namely},
$$
\begin{array}{l}
\nabla_y{\cal L}(x^*,y^*,\mu^*,\lambda^*)=0,\\[3pt]
h(x^*,y^*)=0, \\[3pt]
 0 \leq \lambda^* \perp g(x^*,y^*)\leq 0.
\end{array}
$$
\item[{\rm (b)}] {The linear} independence constraint qualification holds at $y^*$; namely, the set of vectors
$$
\left\{\nabla_y h_1(x^*,y^*),\ldots,\nabla_y h_{m_1}(x^*,y^*)\right\} \cup \left\{\nabla_y g_i(x^*,y^*):i \in I_{x^*}(y^*)\right\}
$$
are linearly independent, where $I_{x^*}(y^*)=\left\{i: g_i(x^*,y^*)=0, i=1,\ldots,m_2\right\}$.
\item[{\rm (c)}] {The strict} complementarity condition holds at $y^*$ for $\lambda^*$; namely,
$$
\lambda^*_i-g_i(x^*,y^*)>0,\quad i=1, \ldots,m_2.
$$
\item[{\rm (d)}] {The second-order} sufficient optimality condition holds at $(y^*,\mu^*,\lambda^*)$,
$$
\langle \nabla_{yy}^2{\cal L}(x^*,y^*,\mu^*,\lambda^*)d_y,d_y\rangle<0 \quad \forall d_y \in {\cal C}_{x^*}(y^*),
$$
where ${\cal C}_{x^*}(y^*)$ is the critical cone of Problem
$({\rm P}_{x^*})$ at $y^*$,
$$
{\cal C}_{x^*}(y^*)=\left\{d_y \in \Re^m: {\cal J}_y h(x^*,y^*)d_y=0; \nabla_yg_i(x^*,y^*)d_y \leq 0, i \in I_{x^*}(y^*);
\nabla_yf(x^*,y^*)d_y\leq 0\right\}.
$$
\end{itemize}
\end{definition}
Let us denote
\begin{equation}\label{def:ab}
\alpha=\Big\{i: g_i(x^*,y^*)=0,i=1,\ldots,p\Big\},\quad \alpha^c=\Big\{i: g_i(x^*,y^*)<0,i=1,\ldots,p\Big\}.
\end{equation}
\begin{lemma}\label{lem:imp}
Let $(x^*,y^*) \in \Re^n \times \Re^m$ be a point around which $f,h,g$ are twice continuously differentiable. Let $(\mu^*,\lambda^*) \in \Re^{m_1}\times \Re^{m_2}$ such that Jacobian uniqueness conditions of Problem
$({\rm P}_{x^*})$ are satisfied at $(y^*,\mu^*,\lambda^*)$. Then there exist $\delta_0>0$ and $\varepsilon_0>0$, and a twice continuously differentiable mapping
$(y,\mu,\lambda):\textbf{B}_{\delta_0}(x^*) \rightarrow \textbf{B}_{\varepsilon_0}(y^*)\times \textbf{B}_{\varepsilon_0}(\mu^*)\times\textbf{B}_{\varepsilon_0}(\lambda^*)$ such that Jacobian uniqueness conditions of Problem $({\rm P}_x)$ are satisfied at $(y(x),\mu(x),\lambda(x))$ when $x \in \textbf{B}_{\delta_0}(x^*)$. Moreover, for  $x\in \textbf{B}_{\delta_0}(x^*)$,
\begin{equation}\label{xa}
\begin{array}{l}
g_i(x,y(x))=0,\ \lambda_i(x)>0,\ i\in \alpha,\\[4pt]
g_i(x,y(x))<0,\  \lambda_i(x)=0,\ i\in \alpha^c.
\end{array}
\end{equation}
\end{lemma}

For $(y(x),\mu(x),\lambda(x))$ given in Lemma \ref{lem:imp}, define the optimal value function
\begin{equation}\label{eq:value-func}
\varphi (x)=f(x,y(x)),\quad x \in \textbf{B}_{\delta_0}(x^*)
\end{equation}
and
\begin{equation}\label{eq:Kx}
K_{\alpha}(x)=
\left
[
\begin{array}{ccc}
\nabla^2_{yy}{\cal L}(x,y(x),\mu(x),\lambda(x)) &  {\cal J}_y h(x,y(x))^T &
   -{\cal J}_y g_{\alpha}(x,y(x))^T\\[4pt]
   {\cal J}_y h(x,y(x)) & 0 &0\\[4pt]
     - {\cal J}_y g_{\alpha}(x,y(x))  & 0& 0
\end{array}
\right
].
\end{equation}

\begin{lemma}\label{ns}
Let $(x^*,y^*) \in \Re^n \times \Re^m$ be a point around which $f,h,g$ are twice continuously differentiable. Let $(\mu^*,\lambda^*) \in \Re^{m_1}\times \Re^{m_2}$ such that Jacobian uniqueness conditions of Problem
$({\rm P}_{x^*})$ are satisfied at $(x^*,\mu^*,\lambda^*)$. Then $K_{\alpha}(x^*)$ is nonsingular and $K_{\alpha}(x)$ is nonsingular when $x \in \textbf{B}_{\delta_0}(x^*)$ for small $\delta_0>0$.
\end{lemma}

Basing on (\ref{xa}), we may  simplify the formula in Proposition 2.1 of \cite{DaiZhang2020} for the second-order derivative of $\varphi(x)$.
\begin{proposition}\label{prop:value-function}
{If} the assumptions of Lemma \ref{lem:imp} are satisfied and $\varphi$ is defined by (\ref{eq:value-func}), then
\begin{equation}\label{phi-1deri}
\nabla_x \varphi(x)=\nabla_x {\cal L}(x,y(x),\mu(x),\lambda(x))
\end{equation}
and
\begin{equation}\label{phi-2deri}
\nabla^2 \varphi(x)=\nabla^2_{xx} {\cal L}(x,y(x),\mu(x),\lambda(x))\\[4pt]
-N_{\alpha}(x)^T K_{\alpha}(x)^{-1}N_{\alpha}(x),
      \end{equation}
where

\begin{equation}\label{Nx}
N_{\alpha}(x)=\left[ \begin{array}{c}
\nabla^2_{x,y}{\cal L}(x,y(x)\mu(x),\lambda(x))\\
   {\cal J}_x h(x,y(x))\\
   {\cal J}_x g_{\alpha}(x,y(x))
   \end{array}
   \right].
\end{equation}
\end{proposition}
{\bf Proof.}  It is easy to check that $K(x)$ is nonsingular when  $x \in \textbf{B}_{\delta_0}(x^*)$ for small $\delta_0>0$. From Proposition 2.1 of  \cite{DaiZhang2020}, we only need to check
\begin{equation}\label{rr}
N(x)^T K(x)^{-1}N(x)=N_{\alpha}(x)^T K_{\alpha}(x)^{-1}N_{\alpha}(x),
\end{equation}
where $K(x)$ and $N(x)$ are defined in \cite{DaiZhang2020} with the following expressions
$$
K(x)=
\left
[
\begin{array}{cccc}
\nabla^2_{yy}{\cal L}(x,y(x),\mu(x),\lambda(x)) & 0 & {\cal J}_y h(x,y(x))^T &
   {\cal J}_y g(x,y(x))^T\\[4pt]
   0  & -2 {\rm Diag}(\lambda(x)) & 0 & 2 {\rm Diag}\left(\sqrt{-g(x,y(x))}\right)\\[4pt]
   {\cal J}_y h(x,y(x)) & 0 &0 &0\\[4pt]
      {\cal J}_y g(x,y(x)) &2 {\rm Diag}\left(\sqrt{-g(x,y(x))}\right) & 0& 0
\end{array}
\right
]
$$
and
$$
N(x)=\left[ \begin{array}{c}
\nabla^2_{x,y}{\cal L}(x,y(x)\mu(x),\lambda(x))\\
 0\\
  {\cal J}_x h(x,y(x))\\
   {\cal J}_x g(x,y(x))
   \end{array}
   \right].
$$
Define
$$
P=\left[
\begin{array}{cccc}
I_{n\times n} &0 &0 &0\\[3pt]
0 & 0 &I_{p\times p} &0\\[3pt]
0 & 0 &0&I_{q\times q}\\[3pt]
0 & I_{p\times p}&0 &0\\[3pt]
\end{array}
\right].
$$
Then $P^TP=PP^T=I_{n+q+2p}$ and
$$
K(x)^{-1}=[P^TPK(x)P^TP]^{-1}=P^T[PK(x)P^T]^{-1}P.
$$
Thus
$$
N(x)^T K(x)^{-1}N(x)=N(x)^TP^T[PK(x)P^T]^{-1}PN(x)=[PN(x)]^T[PK(x)P^T]^{-1}[PN(x)].
$$
Let $G(x)=\nabla^2_{yy}{\cal L}(x,y(x),\mu(x),\lambda(x))$, $D(x)=2 {\rm Diag}\left(\sqrt{-g_{\alpha}(x,y(x))}\right)$ and $E(x)=-2 {\rm Diag}(\lambda_{\alpha}(x))$.
We have that
$$
\begin{array}{l}
PK(x)P^T\\[5pt]
=\left
[
\begin{array}{cccc}
G(x)  & {\cal J}_y h(x,y(x))^T &
   {\cal J}_y g(x,y(x))^T& 0\\[4pt]
   {\cal J}_y h(x,y(x))   & 0 & 0 &0\\[4pt]
  {\cal J}_y g(x,y(x)) & 0 &0 &2 {\rm Diag}\left(\sqrt{-g(x,y(x))}\right)\\[4pt]
     0  &0 & 2 {\rm Diag}\left(\sqrt{-g(x,y(x))}\right)& -2 {\rm Diag}(\lambda(x))
\end{array}
\right
]\\[26pt]
=\left
[
\begin{array}{cccccc}
G(x)  & {\cal J}_y h(x,y(x))^T &
   {\cal J}_y g(x,y(x))^T& 0&0&0\\[4pt]
   {\cal J}_y h(x,y(x))   & 0 & 0 &0&0&0\\[4pt]
  {\cal J}_y g_{\alpha}(x,y(x)) & 0 &0 &0 &0 &0\\[4pt]
  {\cal J}_y g_{\alpha^c}(x,y(x)) & 0 &0 &0&0&D(x)\\[4pt]
  0  &0 & 0 &0& E(x)&0 \\[4pt]
     0  &0 & 0 &D(x)&0 &0
\end{array}
\right
]
\end{array}
$$
Also let
$$
Q=\left[
\begin{array}{cccccc}
I_{n\times} & 0 & 0 &0&0&0\\[3pt]
0 & I_{q\times q} &0&0&0&0\\[3pt]
0& 0& I_{|\alpha|\times |\alpha|} & 0&0&0\\[3pt]
 0& 0& 0&I_{|\alpha^c|\times |\alpha^c|}&0 &0\\[3pt]
0&0&0&0&0 &I_{|\alpha^c|\times |\alpha^c|}\\[3pt]
0&0&0&0&I_{|\alpha|\times |\alpha|} &0
\end{array}
\right].
$$
Then $QQ^T=Q^TQ=I_{n+q+2p}$. Obviously, we have that
$$
\begin{array}{l}
QPK(x)P^TQ^T\\[5pt]
=\left
[
\begin{array}{cccccc}
G(x)  & {\cal J}_y h(x,y(x))^T &
   {\cal J}_y g(x,y(x))^T& 0&0&0\\[4pt]
   {\cal J}_y h(x,y(x))   & 0 & 0 &0&0&0\\[4pt]
  {\cal J}_y g_{\alpha}(x,y(x)) & 0 &0 &0 &0 &0\\[4pt]
  {\cal J}_y g_{\alpha^c}(x,y(x)) & 0 &0 &0&D(x)&0\\[4pt]
  0  &0 & 0 &D(x)& 0&0 \\[4pt]
     0  &0 & 0 &0&0 &E(x)
\end{array}
\right
]
\end{array}
$$
and
$$
QPN(x)=\left[ \begin{array}{c}
\nabla^2_{x,y}{\cal L}(x,y(x)\mu(x),\lambda(x))\\
   {\cal J}_x h(x,y(x))\\
   {\cal J}_x g(x,y(x))\\
   0
   \end{array}
   \right].
$$
Therefore we obtain
$$
\begin{array}{l}
N(x)^T K(x)^{-1}N(x)=[QPN(x)]^T[QPK(x)P^TQ^T]^{-1}[QPN(x)]\\[6pt]
=\left
[
\begin{array}{c}
N_{\alpha}(x)\\[3pt]
{\cal J}_xg_{\alpha^c}(x,y(x))\\[3pt]
0_{\alpha^c}
\end{array}
\right]^T\left[
\begin{array}{ccc}
K_{\alpha}(x) & 0 &0\\[3pt]
0 & 0& D(x)\\[3pt]
0 & D(x)& 0
\end{array}
\right
]^{-1}\left
[
\begin{array}{c}
N_{\alpha}(x)\\[3pt]
{\cal J}_xg_{\alpha^c}(x,y(x))\\[3pt]
0_{\alpha^c}
\end{array}
\right]\\[6pt]
=\left
[
\begin{array}{c}
N_{\alpha}(x)\\[3pt]
{\cal J}_xg_{\alpha^c}(x,y(x))\\[3pt]
0_{\alpha^c}
\end{array}
\right]^T\left[
\begin{array}{ll}
K_{\alpha}(x)^{-1} & 0\quad \quad \quad 0\\[3pt]
\begin{array}{c}
0\\[2pt]
0
\end{array} & \left [
\begin{array}{cc}
0& D(x)\\[3pt]
 D(x)& 0
 \end{array}
 \right
 ]^{-1}
\end{array}
\right
]\left
[
\begin{array}{c}
N_{\alpha}(x)\\[3pt]
{\cal J}_xg_{\alpha^c}(x,y(x))\\[3pt]
0_{\alpha^c}
\end{array}
\right]\\[6pt]
=N_{\alpha}(x)^TK_{\alpha}(x)^{-1}N_{\alpha}(x)+
\left[
\begin{array}{c}
{\cal J}_xg_{\alpha^c}(x,y(x))\\[3pt]
0_{\alpha^c}
\end{array}
\right]^T\left [
\begin{array}{cc}
0& D(x)\\[3pt]
 D(x)& 0
 \end{array}
 \right
 ]^{-1}\left[
\begin{array}{c}
{\cal J}_xg_{\alpha^c}(x,y(x))\\[3pt]
0_{\alpha^c}
\end{array}
\right]\\[8pt]
=N_{\alpha}(x)^TK_{\alpha}(x)^{-1}N_{\alpha}(x).
\end{array}
$$
Namely, (\ref{rr}) holds. The proof is completed.\hfill $\Box$

For $x^* \in \Phi$, the Mangasarian-Fromovitz constraint qualification is said to hold at $x^*$ the constraint set $\Phi$ if
\begin{itemize}
\item[(a)] The set of vectors $\nabla H_j(x^*), j=1,\ldots, n_1$ are linearly independent.
\item[(b)] There exists a vector $\bar d \in \Re^n$ such that
$$
\nabla H_j(x^*)^T\bar d=0, \, j=1,\ldots,n_1,\ \nabla G_i(x^*)^T\bar d<0, \, i \in I(x^*),
$$
where $I(x^*)=\{i: G_i(x^*)=0, i=1,\ldots, n_2\}$.
\end{itemize}
Define the critical cone at $x^*$  by
\begin{equation}\label{eq:critC}
{\cal C}(x^*)=\{d_x \in \Re^n: {\cal J}H(x^*)d_x=0; \nabla G_i(x^*)^Td_x \leq 0, i \in I(x^*); \varphi'(x^*;d_x) \leq0\}.
\end{equation}
In this case, the critical cone ${\cal C}(x^*)$ can be expressed as
\begin{equation}\label{eq:CritiCS}
{\cal C}(x^*)=\{d_x \in \Re^n: {\cal J}H(x^*)d_x=0; \ \nabla G_i(x^*)^Td_x \leq 0, i \in I(x^*); \ \nabla_x {\cal L}(x^*,y^*,\mu^*,\lambda^*)^Td_x\leq 0\}.
\end{equation}
Based on (\ref{rr}) we may simplify Theorems 3.1 and 3.2 in \cite{DaiZhang2020} as follows.
\begin{theorem}\label{th:nc}
(Necessary Optimality Conditions) Let $(x^*,y^*) \in \Re^n \times \Re^m$ be a point around which $f$, $h$, $g$ are twice continuously differentiable and $H$, $G$ are twice continuously differentiable around $x^*$. Let $(x^*,y^*)$  be a local minimax point of Problem (\ref{cminimax}). Assume  that the linear independence constraint qualification holds at $y^*$ for constraint set $Y(x^*)$. Then there exists a unique vector $(\mu^*,\lambda^*) \in \Re^{m_1}\times\Re^{m_2}$ such that
\begin{equation}\label{KKT-Px}
\begin{array}{l}
\nabla_y {\cal L}(x^*,y^*,\mu^*,\lambda^*)=0,\\[3pt]
h(x^*,y^*)=0,\\[3pt]
0\geq \lambda^* \perp g(x^*,y^*) \leq 0.
\end{array}
\end{equation}
For any $d_y \in {\cal C}_{x^*}(y^*)$, we have that
\begin{equation}\label{second-N}
\langle \nabla^2_{yy}{\cal L}(x^*,y^*,\mu^*,\lambda^*)d_y, d_y \rangle \leq 0.
\end{equation}
 Assume further that Problem $({\rm P}_{x^*})$  satisfies  Jacobian uniqueness conditions at $(y^*,\mu^*,\lambda^*)$ and
  the Mangasarian-Fromovitz constraint qualification holds at $x^*$ for the constraint set $\Phi$.
Then there exists $(u^*,v^*) \in \Re^{n_1}\times \Re^{n_2}$ such that
 \begin{equation}\label{KKT-P}
\begin{array}{l}
\nabla_x {\cal L}(x^*,y^*,\mu^*,\lambda^*)+{\cal J}H(x^*)^Tu^*+{\cal J}G(x^*)^Tv^*=0,\\[3pt]
H(x^*)=0,\\[3pt]
0\leq v^* \perp G(x^*) \leq 0.
\end{array}
\end{equation}
The set of all $(u^*,v^*)$ satisfying (\ref{KKT-P}), denoted by $\Lambda (x^*)$, is nonempty compact convex set.
Furthermore, for every $d_x \in {\cal C}(x^*)$, where ${\cal C}(x^*)$ is defined by (\ref{eq:CritiCS}),
\begin{equation}\label{eq:secNCs}
\begin{array}{l}
\displaystyle \max_{(u,v) \in \Lambda(x^*)} \left\{ \left\langle \left[\displaystyle \sum_{j=1}^{n_1}u_i\nabla^2_{xx}H_j(x^*)+\displaystyle \sum_{i=1}^{n_2} v_i\nabla^2_{xx}G_i(x^*)\right]d_x,d_x \right\rangle\right\}\\[16pt]
 \quad \quad \quad \, +\left\langle \left[ \nabla^2_{xx}{\cal L}(x^*,y^*,\mu^*,\lambda^*)- N_{\alpha}(x^*)^TK_{\alpha}(x^*)^{-1}N_{\alpha}(x^*)\right]d_x,d_x \right\rangle\geq0,
 \end{array}
\end{equation}
where $K_{\alpha}(x)$ is defined by (\ref{eq:Kx}) and $N_{\alpha}(x)$ is defined by (\ref{Nx}).

\end{theorem}

We name the first-order necessary optimality conditions (\ref{KKT-P}) and (\ref{KKT-Px})
as KKT conditions of Problem (\ref{cminimax}) at $(x^*,u^*,v^*,y^*,\mu^*,\lambda^*)$.

\begin{theorem}\label{th:Sc}
(Second-order Sufficient Optimality Conditions)
Let $(x^*,y^*) \in \Re^n \times \Re^m$ be a point around which $f,h,g$ are twice continuously differentiable and $H$, $G$ are twice continuously differentiable around $x^*$. Assume that $x^* \in \Phi$ and $y^* \in Y(x^*)$.  Let $(\mu^*,\lambda^*) \in \Re^{m_1}\times \Re^{m_2}$.
 Suppose that  Problem $({\rm P}_{x^*})$  satisfies  Jacobian uniqueness conditions at $(y^*,\mu^*,\lambda^*)$, $\Lambda (x^*) \ne \emptyset$, and for every $d_x \in {\cal C}(x^*)\setminus \emptyset$ (where ${\cal C}(x^*)$ is defined by (\ref{eq:CritiCS})),
\begin{equation}\label{eq:secSCs}
\begin{array}{l}
\displaystyle \sup_{(u,v) \in \Lambda(x^*)} \left\{ \left\langle \left[\displaystyle \sum_{j=1}^{n_1}u_i\nabla^2_{xx}H_j(x^*)+\displaystyle \sum_{i=1}^{n_2} v_i\nabla^2_{xx}G_i(x^*)\right]d_x,d_x \right\rangle\right\}\\[16pt]
 \quad \quad \quad \, +\left\langle \left[ \nabla^2_{xx}{\cal L}(x^*,y^*,\mu^*,\lambda^*)- N_{\alpha}(x^*)^TK_{\alpha}(x^*)^{-1}N_{\alpha}(x^*)\right]d_x,d_x \right\rangle > 0,
 \end{array}
\end{equation}
where $K_{\alpha}(x)$ is defined by (\ref{eq:Kx}) and $N_{\alpha}(x)$ is defined by (\ref{Nx}).
Then there exist $\delta_1 \in (0,\delta_0)$, $\varepsilon_1 \in (0,\varepsilon_0)$ (where $\delta_0$ and $\varepsilon_0$ are given by Lemma \ref{lem:imp}) and $\gamma_1>0$,$\gamma_2>0$ such that
for $x \in \textbf{B}_{\delta_1}(x^*)\cap \Phi$ and $y \in \textbf{B}_{\varepsilon_1}(y^*)\cap Y(x^*)$,
\begin{equation}\label{eq:2ndG}
f(x^*,y)+\gamma_1 \|y-y^*\|^2/2 \leq f(x^*,y^*)\leq \displaystyle \sup_{z \in Y(x) \cap \textbf{B}_{\varepsilon_0}(y^*)} f(x,z)-\gamma_2\|x-x^*\|^2/2,
\end{equation}
which implies that $(x^*,y^*)$ is a local minimax point of Problem (\ref{cminimax}).
\end{theorem}

\section{Stability under Jacobian Uniqueness Condition}
\setcounter{equation}{0}
\quad \, For convenience in stating the  stability result about ${\cal C}^2$ perturbation of Problem (\ref{cminimax}) when the conditions in Theorem \ref{th:Sc} are satisfied, we introduce the following definition.
\begin{definition}\label{def:CS}
If the following conditions are satisfied, we say that Problem (\ref{cminimax}) satisfies  Jacobian uniqueness condition at $(x^*,u^*,v^*,y^*,\mu^*,\lambda^*)\in \Re^n\times \Re^{n_1}\times \Re^{n_2}\times \Re^m \times \Re^{m_1}\times \Re^{m_2}$.
\begin{itemize}
\item[{\rm (i)}]$x^* \in \Phi$ and conditions in (\ref{KKT-P}) are satisfied at $(x^*,u^*,v^*,y^*,\mu^*,\lambda^*)$.
\item[{\rm (ii)}]The set vectors $\big\{\nabla H_1(x^*),\ldots,\nabla H_{n_1}\big\}\cup \Big\{\nabla G_i(x^*):i\in I(x^*)\Big\}$
are linearly independent, where $I(x^*)=\{i: G_i(x^*)=0, i=1.\ldots, n_2\}$.
\item[{\rm (iii)}]$v^*_i-G_i(x^*)>0$ for $i \in I(x^*)$.
\item[{\rm (iv)}]$y^* \in Y(x^*)$ and  Problem $({\rm P}_{x^*})$  satisfies  Jacobian uniqueness conditions at $(y^*,\mu^*,\lambda^*)$.
    \item[{\rm (v)}]For every $d_x \in {\cal C}(x^*)\setminus \{0\}$ (where ${\cal C}(x^*)$ is defined by (\ref{eq:CritiCS})), the second-order sufficient optimality condition (\ref{eq:secSCs}) is satisfied.
\end{itemize}
\end{definition}
Following Kojima (1980)
\cite{Kojima80}, we define the so-called Kojima  mapping for  Problem (\ref{cminimax}),
\begin{equation}\label{Kojima-f}
F(x,u,w,y,\mu,\xi)=\left[
\begin{array}{c}
\nabla_x{\cal L}(x,y,\mu,\xi^+)+{\cal J}H(x)^Tu+{\cal J}G(x)^Tw^+\\[3pt]
H(x)\\[3pt]
G(x)-w^-\\[3pt]
\nabla_y{\cal L}(x,y,\mu,\xi^+)\\[3pt]
h(x,y)\\[3pt]
-g(x,y)+\xi^-
\end{array}
\right
],
\end{equation}
where
 $w^+_i=\max\,\big\{0,w_i\big\}$, $\xi^-_i=\min\,\big\{0,w_i\big\}$, $i=1,\ldots, n_2$ for $w \in \Re^{n_2}$ and
 $\xi^+_i=\max\,\big\{0,\xi_i\big\}$, $\xi^-_i=\min\,\big\{0,\xi_i\big\}$, $i=1,\ldots,m_2$ for $\xi \in \Re^{m_2}$.
 If $F(x,u,w,y,\mu,\xi)=0$, then, letting $v=w^+$ and $\lambda=\xi^+$, $(x,y,u,v,\mu,\lambda)$ satisfies
 the first-order necessary optimality conditions of Problem (\ref{cminimax}).

\begin{lemma}\label{lem:nsing}
Suppose that the Jacobian uniqueness condition of Problem (\ref{cminimax}) in Definition \ref{def:CS} is satisfied at $(x^*,u^*,v^*,y^*,\mu^*,\lambda^*)\in \Re^n\times \Re^{n_1}\times \Re^{n_2}\times \Re^m \times \Re^{m_1}\times \Re^{m_2}$. Then $F$ is differentiable at $(x^*,u^*,w^*,y^*,\mu^*,\xi^*)$ and ${\cal J}F(x^*,u^*,w^*,y^*,\mu^*,\xi^*)$ is nonsingular for $w^*=v^*+G(x^*)$ and
$\xi^*=\lambda^*+g(x^*,y^*)$.
\end{lemma}
{\bf Proof}. Since $\lambda^*-g(x^*,y^*)>0$ and $v^*-G(x^*)>0$, we know that $\xi^+$ and $\xi^-$ are differentiable at $\xi^*$, and $w^+$ and $w^-$ are differentiable at $w^*$. Thus $F$ is differentiable at $(x^*,y^*,u^*,w^*,\mu^*,\xi^*)$. Without loss of generality, we assume that
$$
\beta:=I(x^*)=\{1,\ldots, r\},\quad \alpha=\{1,\ldots, s\},
$$
where
$$
I(x^*)=\{i:G_i(x^*)=0, i=1,\ldots, n_2\},\quad \alpha=\{i: g_i(x^*,y^*)=0,i=1,\ldots,m_2\}.
$$
Then, for $\beta^c=\{1,\ldots,n_2\}\setminus \beta$ and $\alpha^c=\{1,\ldots,m_2\}\setminus \alpha$, we get that
$$
\beta^c=\{r+1,\ldots, n_2\}, \quad \alpha^c=\{s+1,\ldots, m_2\}.
$$
Thus we obtain
\begin{equation}\label{eq:J}
\begin{array}{cclccl}
{\cal J}w^+|_{w=w^*}&=&\left[
\begin{array}{cc}
I_r & 0\\[3pt]
0 & 0
\end{array}
\right],\quad & {\cal J}w^-|_{w=w^*}&=&\left[
\begin{array}{cc}
0 & 0\\[3pt]
0 & I_{n_2-r}
\end{array}
\right],\\[8pt]
{\cal J}\xi^+|_{\xi=\xi^*}&=&\left[
\begin{array}{cc}
I_s & 0\\[3pt]
0 & 0
\end{array}
\right],\quad & {\cal J}\xi^-|_{\xi=\xi^*}&=&\left[
\begin{array}{cc}
0 & 0\\[3pt]
0 & I_{m_2-s}
\end{array}
\right].
\end{array}
\end{equation}
Denote
$$
\begin{array}{ll}
G^*_{11}=\nabla^2_{xx}{\cal L}(x^*,y^*,\mu^*,\lambda^*)
+\displaystyle \sum_{j=1}^{n_1}u_i\nabla^2_{xx}H_j(x^*)+\displaystyle \sum_{i=1}^{n_2} v_i\nabla^2_{xx}G_i(x^*),\\[4pt]
G^*_{12}=\nabla^2_{xy}{\cal L}(x^*,y^*,\mu^*,\lambda^*),\\[4pt]
G^*_{22}=\nabla^2_{yy}{\cal L}(x^*,y^*,\mu^*,\lambda^*).
 \end{array}
$$
For simplicity, we use notations ${\cal J}_xh^*$ and ${\cal J}_yh^*$ to represent ${\cal J}_xh(x^*,y^*)$ and ${\cal J}_yh(x^*,y^*)$, respectively. The same
notations are also applied to $g_{\alpha}$ and $g_{\alpha^c}$. Then the Jacobian of $F$ at $(x^*,u^*,w^*,y^*,\mu^*,\xi^*)$ can be expressed as
\begin{equation}\label{JF}
\begin{array}{l}
{\cal J}F(x^*,u^*,w^*,y^*,\mu^*,\xi^*)\\[6pt]
=\left[
\begin{array}{cccccccc}
G^*_{11}& {\cal J}H(x^*)^T & {\cal J}G_{\beta}(x^*)^T &0 &G^*_{12}&{\cal J}_xh^{*T}&-{\cal J}_xg_{\alpha}^{*T} &0\\[6pt]
{\cal J}H(x^*) & 0 &0 &0 &0 &0 &0 &0 \\[6pt]
{\cal J}G_{\beta}(x^*) & 0&0&0&0&0&0&0\\[6pt]
{\cal J}G_{\beta^c}(x^*) & 0&0&-I_{n_2-r}&0&0&0&0\\[6pt]
G^{*T}_{12} & 0&0&0&G^*_{22}&{\cal J}_yh^{*T}&-{\cal J}_yg_{\alpha}^{*T} &0\\[6pt]
{\cal J}_xh^* & 0&0&0&{\cal J}_yh^*&0&0&0\\[6pt]
-{\cal J}_xg_{\alpha}^* & 0&0&0&-{\cal J}_yg_{\alpha}^*&0&0&0\\[6pt]
-{\cal J}_xg_{\alpha^c}^* & 0&0&0&-{\cal J}_yg_{\alpha^c}^*&0&0&I_{m_2-s}
\end{array}
\right].
\end{array}
\end{equation}
The nonsingularity of ${\cal J}F(x^*,u^*,w^*,y^*,\mu^*,\xi^*)$ is equivalent to the nonsingularity of the following matrix
$$
\begin{array}{l}
\left[
\begin{array}{cccccccc}
G^*_{11}& {\cal J}H(x^*)^T & {\cal J}G_{\beta}(x^*)^T &G^*_{12}&{\cal J}_xh^{*T}&-{\cal J}_xg_{\alpha}^{*T} &0&0 \\[6pt]
{\cal J}H(x^*) & 0 &0 &0 &0 &0 &0 &0 \\[6pt]
{\cal J}G_{\beta}(x^*) & 0&0&0&0&0&0&0\\[6pt]
G^{*T}_{12} & 0&0&G^*_{22}&{\cal J}_yh^{*T}&-{\cal J}_yg_{\alpha}^{*T} &0&0\\[6pt]
{\cal J}_xh^* & 0&0&{\cal J}_yh^*&0&0&0&0\\[6pt]
-{\cal J}_xg_{\alpha}^* & 0&0&-{\cal J}_yg_{\alpha}^*&0&0&0&0\\[6pt]
{\cal J}G_{\beta^c}(x^*) & 0&0&0&0&0&-I_{n_2-r}&0\\[6pt]
-{\cal J}_xg_{\alpha^c}^* & 0&0&-{\cal J}_yg_{\alpha^c}^*&0&0&0&I_{m_2-s}
\end{array}
\right],
\end{array}
$$
which is equivalent to the nonsingularity of the following
matrix
\begin{equation}\label{JFa}
H=
\left[
\begin{array}{cccccc}
G^*_{11}& {\cal J}H(x^*)^T & {\cal J}G_{\beta}(x^*)^T &G^*_{12}&{\cal J}_xh^{*T}&-{\cal J}_xg_{\alpha}^{*T}  \\[6pt]
{\cal J}H(x^*) & 0 &0 &0&0&0  \\[6pt]
{\cal J}G_{\beta}(x^*) & 0&0&0&0&0\\[6pt]
G^{*T}_{12} & 0&0&G^*_{22}&{\cal J}_yh^{*T}&-{\cal J}_yg_{\alpha}^{*T}\\[6pt]
{\cal J}_xh^* & 0&0&{\cal J}_yh^*&0&0\\[6pt]
-{\cal J}_xg_{\alpha}^* & 0&0&-{\cal J}_yg_{\alpha}^*&0&0
\end{array}
\right].
\end{equation}
Therefore, we only need to prove that the matrix $H$ is nonsingular. From Lemma \ref{ns}, we obtain
that $K_{\alpha}(x^*)$ is nonsingular, where
$$
K_{\alpha}(x^*)=\left[
\begin{array}{ccc}
G^*_{22}&{\cal J}_yh^{*T}&-{\cal J}_yg_{\alpha}^{*T}\\[6pt]
{\cal J}_yh^*&0&0\\[6pt]
-{\cal J}_yg_{\alpha}^*&0&0
\end{array}
\right].
$$
So it suffices to prove that
$
H/K_{\alpha}(x^*)
$ is nonsingular. Noticing that
\begin{equation}\label{Hkf}
\begin{array}{rcl}
H/K_{\alpha}(x^*)&=& \left[
\begin{array}{ccc}
G^*_{11}& {\cal J}H(x^*)^T & {\cal J}G_{\beta}(x^*)^T  \\[6pt]
{\cal J}H(x^*) & 0 &0 \\[6pt]
{\cal J}G_{\beta}(x^*) & 0&0
\end{array}
\right]\\[20pt]
&&-\left[
\begin{array}{ccc}
G^*_{12}&{\cal J}_xh^{*T}&-{\cal J}_xg_{\alpha}^{*T}  \\[6pt]
0 &0&0  \\[6pt]
0&0&0
\end{array}
\right]K_{\alpha}(x^*)^{-1}\left[
\begin{array}{ccc}
G^{*T}_{12}&0&0  \\[6pt]
{\cal J}_xh^* &0&0  \\[6pt]
-{\cal J}_xg_{\alpha}^*&0&0
\end{array}
\right]\\[20pt]
&=&\left[
\begin{array}{ccc}
G^*_{11}-N_{\alpha}(x^*)^TK_{\alpha}(x^*)N_{\alpha}(x^*)& {\cal J}H(x^*)^T & {\cal J}G_{\beta}(x^*)^T  \\[6pt]
{\cal J}H(x^*) & 0 &0 \\[6pt]
{\cal J}G_{\beta}(x^*) & 0&0
\end{array}
\right],
\end{array}
\end{equation}
we have from (iii) that ${\cal C}(x^*)$ is reduced to the following subspace
\begin{equation}\label{eq:kerhg}
{\cal C}(x^*)={\rm Ker}\, \left[
\begin{array}{l}
{\cal J}H(x^*) \\[6pt]
{\cal J}G_{\beta}(x^*)
\end{array}
\right
].
\end{equation}
Now we prove that $H/K_{\alpha}(x^*)$ is nonsingular via the formula (\ref{Hkf}). Let $a\in \Re^n$, $b\in \Re^{n_1}$ and $
c\in \Re^r$ satisfy
$$
\left[
\begin{array}{ccc}
G^*_{11}-N_{\alpha}(x^*)^TK_{\alpha}(x^*)N_{\alpha}(x^*)& {\cal J}H(x^*)^T & {\cal J}G_{\beta}(x^*)^T  \\[6pt]
{\cal J}H(x^*) & 0 &0 \\[6pt]
{\cal J}G_{\beta}(x^*) & 0&0
\end{array}
\right]\left[
\begin{array}{c}
a\\[6pt]
b\\[6pt]
c
\end{array}
\right
]=0
$$
or
\begin{equation}\label{abc}
\left\{
\begin{array}{l}
G^*_{11}-N_{\alpha}(x^*)^TK_{\alpha}(x^*)N_{\alpha}(x^*)a+ {\cal J}H(x^*)^T b+ {\cal J}G_{\beta}(x^*)^T c=0, \\[6pt]
{\cal J}H(x^*)a=0, \\[6pt]
{\cal J}G_{\beta}(x^*)a=0.
\end{array}
\right.
\end{equation}
It follows from ${\cal J}H(x^*)a=0$ and ${\cal J}G_{\beta}(x^*) a=0$ that $a \in {\cal C}(x^*)$.
Premultiplying $a^T$ to the first equation in (\ref{abc}), we obtain
$$
\begin{array}{l}
\displaystyle  \left\langle \left[\displaystyle \sum_{j=1}^{n_1}u^*_i\nabla^2_{xx}H_j(x^*)+\displaystyle \sum_{i=1}^{n_2} v^*_i\nabla^2_{xx}G_i(x^*)\right]a,a \right\rangle\\[16pt]
 \quad \quad \quad \, +\left\langle \left[ \nabla^2_{xx}{\cal L}(x^*,y^*,\mu^*,\lambda^*)- N_{\alpha}(x^*)^TK_{\alpha}(x^*)^{-1}N_{\alpha}(x^*)\right]a,a \right\rangle= 0,
 \end{array}
$$
which implies $a=0$ from the condition (\ref{eq:secSCs}). From the first equation in (\ref{abc}) again, we obtain
$$
{\cal J}H(x^*)^T b+ {\cal J}G_{\beta}(x^*)^T c=0,
$$
from which we obtain $b=0$ and $c=0$ from
 (iii).  Therefore $H/K_{\alpha}(x^*)$ is nonsingular.  The proof is completed. \hfill $\Box$

Basing on Lemma \ref{lem:nsing}, we may establish the stability on the ${\cal C}^2$ perturbation of Problem (\ref{cminimax})
under the Jacobian uniqueness condition  at $(x^*,u^*,v^*,y^*,\mu^*,\lambda^*)$.

Now consider the parameterized  constrained minimax optimization problem of the form
\begin{equation}\label{Pcminimax}
({\rm P}_\vartheta) \quad \quad \quad \min_{x \in \Phi(\vartheta)}\max_{y \in \bar Y(x,\vartheta)}\bar f(x,y,\vartheta),
\end{equation}
where $\bar f:\Re^n\times \Re^m \times \Re^l\rightarrow \Re$, $\Phi\subset \Re^n$ is a feasible set of decision variable $x$ defined by
\begin{equation}\label{PPhi}
\Phi (\vartheta)=\{x \in \Re^n: \bar H(x,\vartheta)=0,\, \bar G(x,\vartheta)\leq 0\}
\end{equation}
and $Y: \Re^n \times \Re^l \rightrightarrows \Re^m$ is a set-valued mapping defined by
\begin{equation}\label{PYx}
\bar Y(x,\vartheta)=\{y \in \Re^m: \bar h(x,y,\vartheta)=0,\, \bar g(x,y,\vartheta)\leq 0\}.
\end{equation}
Let $\vartheta_0\in \Re^l$ be a vector such that
$$
f(x,y)=\bar f(x,y,\vartheta_0),\quad h(x,y)=\bar h(x,y,\vartheta_0),\quad g(x,y)=\bar g(x,y,\vartheta_0)
$$
and
$$
H(x)=\bar H(x,\vartheta_0),\quad G(x)=\bar G(x,\vartheta_0).
$$
\begin{definition}\label{def:c2P}
We say Problem (\ref{Pcminimax}) is a local ${\cal C}^2$ perturbation of Problem (\ref{cminimax}) around $(x^*,y^*)$ if there exist  open sets ${\cal O}_1\subset \Re^n$, ${\cal O}_2\subset \Re^m$ and  $\Theta\subset \Re^l$ satisfying $\vartheta_0\in \Theta$, $x^* \in {\cal O}_1$, $y^* \in {\cal O}_2$ and
$\bar f$, $\bar h$, $\bar g$ are twicely smooth over ${\cal O}_1\times {\cal O}_2\times \Theta$, and $\bar H$, $\bar G$
 are twicely smooth over ${\cal O}_1\times \Theta$.
\end{definition}
The Kojima  mapping for  Problem $({\rm {P}}_{\vartheta})$ is the following function
\begin{equation}\label{PKojima-f}
\bar F(x,u,w,y,\mu,\xi;\vartheta)=\left[
\begin{array}{c}
\nabla_x\bar{\cal L}(x,y,\mu,\xi^+;\vartheta)+{\cal J}H(x,\vartheta)^Tu+{\cal J}G(x,\vartheta)^Tw^+\\[3pt]
H(x,\vartheta)\\[3pt]
G(x,\vartheta)-w^-\\[3pt]
\nabla_y\bar{\cal L}(x,y,\mu,\xi^+;\vartheta)\\[3pt]
h(x,y,\vartheta)\\[3pt]
-g(x,y,\vartheta)+\xi^-
\end{array}
\right
],
\end{equation}
where
$$
\bar {\cal L}(x,y,\mu,\lambda;\vartheta)=\bar f(x,y,\vartheta)+\langle  \mu, \bar h(x,y,\vartheta) \rangle
-\langle  \lambda, \bar g(x,y,\vartheta) \rangle.
$$
\begin{theorem}\label{th1}
Suppose that the Jacobian uniqueness condition of Problem (\ref{cminimax}) is satisfied at $(x^*$, $u^*$, $v^*$, $y^*$, $\mu^*$, $\lambda^*)\in \Re^n\times \Re^{n_1}\times \Re^{n_2}\times \Re^m \times \Re^{m_1}\times \Re^{m_2}$ and (${\rm P}_{\vartheta}$) is a local ${\cal C}^2$ perturbation of Problem (\ref{cminimax}) around $(x^*,y^*)$.
 Then there exist $\varepsilon>0$ and $\delta>0$ such that $\textbf{B}(\vartheta_0,\delta)\subset \Theta$,$\textbf{B}(x^*,\varepsilon)\subset {\cal O}_1$ and $\textbf{B}(y^*,\varepsilon)\subset {\cal O}_2$, and there is a mapping $(x(\cdot),y(\cdot),u(\cdot),v(\cdot),\mu(\cdot),\lambda(\cdot)): \textbf{B}(\vartheta_0,\delta)\rightarrow
 \textbf{B}(x^*,\varepsilon) \times \textbf{B}(y^*,\varepsilon)\times \textbf{B}(u^*,\varepsilon)\times \textbf{B}(v^*,\varepsilon)\times \textbf{B}(\mu^*,\varepsilon)\times \textbf{B}(\lambda^*,\varepsilon)$ such that, for
 $v(\vartheta)=w(\vartheta)^+$ and $\lambda(\vartheta)=\xi(\vartheta)^+$,
 \begin{itemize}
 \item[{\rm (1)}]$(x(\vartheta_0),y(\vartheta_0),u(\vartheta_0),v(\vartheta_0),\mu(\vartheta_0),\lambda(\vartheta_0))=
 (x^*,y^*,u^*,v^*,\mu^*,\lambda^*)$.
 \item[{\rm (2)}]For any $\vartheta \in  \textbf{B}(\vartheta_0,\delta)$, $(x(\cdot),y(\cdot),u(\cdot),v(\cdot),\mu(\cdot),\lambda(\cdot))$ is continuously differentiable at $\vartheta$.
 \item[{\rm (3)}]For any $\vartheta \in  \textbf{B}(\vartheta_0,\delta)$, Problem (${\rm P}_{\vartheta}$)  satisfies
 the Jacobian uniqueness condition at $(x(\vartheta),y(\vartheta),$ $u(\vartheta),v(\vartheta),\mu(\vartheta),\lambda(\vartheta))$.
 \end{itemize}
  \end{theorem}
  {\bf Proof}. Let $w^*=v^*+G(x^*)$ and
$\xi^*=\lambda^*+g(x^*,y^*)$. From the definitions of $\bar F$ and $F$ in (\ref{Kojima-f}), we have
$$
\bar F(x,u,w,y,\mu,\xi;\vartheta_0)=F(x,u,w,y,\mu,\xi).
$$
Thus we get that
$$\bar F(x^*,u^*,w^*,y^*,\mu^*,\xi^*;\vartheta_0)=0,\ {\cal J}_{(x,u,w,y,\mu,\xi)}\bar F(x^*,u^*,w^*,y^*,\mu^*,\xi^*;\vartheta_0)={\cal J}F(x^*,u^*,w^*,y^*,\mu^*,\xi^*),$$
 and in turn  ${\cal J}_{(x,u,w,y,\mu,\xi)}\bar F(x^*,u^*,w^*,y^*,\mu^*,\xi^*;\vartheta_0)$ is nonsingular from Lemma \ref{lem:nsing}. From the classical implicit function theorem, we get that there exist $\varepsilon>0$ and $\delta>0$ such that $\textbf{B}(\vartheta_0,\delta)\subset \Theta$,$\textbf{B}(x^*,\varepsilon)\subset {\cal O}_1$ and $\textbf{B}(y^*,\varepsilon)\subset {\cal O}_2$, and there is a mapping $(x(\cdot),u(\cdot),w(\cdot),y(\cdot),\mu(\cdot),\xi(\cdot)): \textbf{B}(\vartheta_0,\delta)\rightarrow
 \textbf{B}(x^*,\varepsilon) \times \textbf{B}(u^*,\varepsilon)\times \textbf{B}(w^*,\varepsilon)\times \textbf{B}(y^*,\varepsilon)\times \textbf{B}(\mu^*,\varepsilon)\times \textbf{B}(\xi^*,\varepsilon)$ such that
$$
(x(\vartheta_0),u(\vartheta_0),w(\vartheta_0),y(\vartheta_0),\mu(\vartheta_0),\xi(\vartheta_0))=
 (x^*,u^*,w^*,y^*,\mu^*,\xi^*);
$$
meanwhile, for any $\vartheta \in  \textbf{B}(\vartheta_0,\delta)$, $(x(\cdot),u(\cdot),w(\cdot),y(\cdot),\mu(\cdot),\xi(\cdot))$ is continuously differentiable at $\vartheta$, and
\begin{equation}\label{eqtheta}
\bar F(x(\vartheta),u(\vartheta),w(\vartheta),y(\vartheta),\mu(\vartheta),\xi(\vartheta);\vartheta)=0,\,\, \forall \vartheta \in \textbf{B}(\vartheta_0,\delta).
\end{equation}
From the continuity of $(x(\vartheta),u(\vartheta),w(\vartheta),y(\vartheta),\mu(\vartheta),\xi(\vartheta))$ for $\vartheta \in \textbf{B}(\vartheta_0,\delta)$, we have from $G_i(x(\vartheta))-w_i(\vartheta)^-=0$ in (\ref{eqtheta}) that
$$
G_i(x(\vartheta))<0, \quad w^+(\vartheta)_i=0, i \in \beta^c
$$
 and
 $$
 w^+_i(\vartheta)>0, \quad G_i(x(\vartheta))=0, i \in \beta.
 $$
 Therefore we have that
 $$
\beta(\vartheta):=\{i: G_i(x(\vartheta))=0, i=1,\ldots,n_2\}=\beta, \quad \{i: w^+_i(x(\vartheta))=0, i=1,\ldots,n_2\}=\beta^c
 $$
 and that $w^+_i(\cdot)$ is differentiable at $\vartheta$ for $i\in \beta$ and $w^+_i(\cdot)\equiv0 $ for $i \in \beta^c$.
In turn $v(\vartheta)=w^+(\vartheta)$ is differentiable over $\textbf{B}(\vartheta_0,\delta)$.
Using the same arguments as the above, we  obtain
   $$
\alpha(\vartheta):=\{i: g_i(x(\vartheta),y(\vartheta))=0, i=1,\ldots,m_2\}=\alpha, \quad \{i: \xi^+_i(x(\vartheta))=0, i=1,\ldots,m_2\}=\alpha^c
 $$
 and
 $\lambda (\vartheta)=\xi(\vartheta)^+$ is also differentiable over $\textbf{B}(\vartheta_0,\delta)$.  Hence the assertions (1) and (2) hold.

  Now we prove the assertion (3). From the first three equations in (\ref{eqtheta}) and the definition of $(v(\vartheta),\lambda(\vartheta))$, we obtain
$$
\begin{array}{l}
\nabla_x{\cal L}(x(\vartheta),y(\vartheta),\mu(\vartheta),\lambda(\vartheta))+{\cal J}H(x(\vartheta))^Tu(\vartheta)+{\cal J}G(x(\vartheta))^Tv(\vartheta)=0,\\[3pt]
H(x(\vartheta))=0,\\[3pt]
G(x(\vartheta))-[G(x(\vartheta)+v(\vartheta)]^-=0, \\[3pt]
 \end{array}
 $$
 which are exactly the conditions in  (i) of the Jacobian uniqueness condition of $({\rm P}_{\vartheta})$ from Definition
 \ref{def:CS}. From the continuity of $(x(\vartheta),v(\vartheta))$ and $\beta(\vartheta)=\beta$ claimed just now, we have that
 the set vectors $\big\{\nabla H_1(x(\vartheta)),\ldots,\nabla H_{n_1}(x(\vartheta))\big\}\cup \Big\{\nabla G_i(x(\vartheta)):i\in \beta(\vartheta)\Big\}$
are linearly independent, where $\beta(\vartheta)=\{i: G_i(x(\vartheta))=0, i=1.\ldots, n_2\}$, which implies  (ii) of  the Jacobian uniqueness condition of $({\rm P}_{\vartheta})$.  And we have that $v_i(\vartheta)-G_i(x(\vartheta))>0$ for $i \in \beta (\vartheta)$; namely, (iii) of  the Jacobian uniqueness condition of $({\rm P}_{\vartheta})$ holds.

  Now we check (iv)  of  the Jacobian uniqueness condition of $({\rm P}_{x(\vartheta)})$; namely,  (${\rm P}_{x(\vartheta)}$) satisfies the Jacobian uniqueness condition at $(y(\vartheta),\mu(\vartheta),\lambda(\vartheta))$. From the last three equations in (\ref{eqtheta}) and the definition of $(v(\vartheta),\lambda(\vartheta))$, we obtain
$$
\begin{array}{l}
\nabla_y{\cal L}(x(\vartheta),y(\vartheta),\mu(\vartheta),\lambda(\vartheta))=0,\\[3pt]
h(x(\vartheta),y(\vartheta))=0,\\[3pt]
g(x(\vartheta),y(\vartheta))-[g(x(\vartheta),y(\vartheta))+\lambda(\vartheta)]^-=0,\\[3pt]
 \end{array}
 $$
which are just KKT conditions of  $({\rm P}_{x(\vartheta)})$ at $(y(\vartheta)$, $\mu(\vartheta)$, $\lambda(\vartheta))$. Since the continuity of $(x(\vartheta)$, $y(\vartheta)$, $\mu(\vartheta)$, $\lambda(\vartheta))$ and $\alpha(\vartheta)=\alpha$, we have that
 the set vectors $$\big\{\nabla_y h_1(x(\vartheta),y(\vartheta)),\ldots,\nabla_y h_{n_1}(x(\vartheta),y(\vartheta))\big\}\cup \Big\{\nabla_y g_i(x(\vartheta),y(\vartheta)):i\in \alpha(\vartheta)\Big\}$$ are linearly independent when $\delta>0$ is small enough; namely, the linear independence constraint qualification of $({\rm P}_{x(\vartheta)})$ at $y(\vartheta)$ is satisfied. We also have $\lambda_i(\vartheta)-g_i(x(\theta),y(\vartheta))>0$; namely, the strict complementarity condition of
 $({\rm P}_{x(\vartheta)})$ holds at $(y(\vartheta),\lambda(\vartheta))$. Until now, for the Jacobian uniqueness condition of $({\rm P}_{x(\vartheta)})$, only the second-order sufficient optimality condition is left to prove.  It can be proved in the same way as that for (v) of the Jacobian uniqueness condition of $({\rm P}_{\vartheta})$ from Definition
 \ref{def:CS}. We omit it here.

 Finally, we prove (v) of the Jacobian uniqueness condition of $({\rm P}_{\vartheta})$  at $(x(\vartheta)$, $y(\vartheta)$, $u(\vartheta)$, $v(\vartheta)$, $\mu(\vartheta)$, $\lambda(\vartheta))$ from Definition
 \ref{def:CS}. From $\alpha(\vartheta)=\alpha$, we have that
 \begin{equation}\label{cxtheta}
 {\cal C}_{x(\vartheta)}(y(\vartheta))={\rm ker}\, {\cal J}_yh(x(\vartheta),y(\vartheta))\cap {\rm ker}\, {\cal J}_yg_{\alpha}(x(\vartheta),y(\vartheta))
 \end{equation}
is a subspace of $\Re^m$. For proving  (v) of the Jacobian uniqueness condition of $({\rm P}_{\vartheta})$  at $(x(\vartheta),y(\vartheta),$ $u(\vartheta),v(\vartheta),\mu(\vartheta),\lambda(\vartheta))$ from Definition
 \ref{def:CS}, we only need to construct a matrix $Z(\vartheta) \in \Re^{m,m-m_1-|\alpha|}$ such that
 \begin{equation}\label{ta1}
 {\rm Range}\, Z(\vartheta)={\rm ker}\, {\cal J}_yh(x(\vartheta),y(\vartheta))\cap {\rm ker}\, {\cal J}_yg_{\alpha}(x(\vartheta),y(\vartheta))
 \end{equation}
 and
  \begin{equation}\label{ta2}
 Z(\vartheta)^T\Psi(\vartheta)Z(\vartheta) \succ 0,
 \end{equation}
 where
 $$
\begin{array}{ll}
\Psi(\vartheta)&=\displaystyle \sum_{j=1}^{n_1}u_i(\vartheta)\nabla^2_{xx}H_j(x(\vartheta))+\displaystyle \sum_{i=1}^{n_2} v_i(\vartheta)\nabla^2_{xx}G_i(x(\vartheta))\\[16pt]
 &\quad \quad \quad \, + \nabla^2_{xx}{\cal L}(x(\vartheta),y(\vartheta),\mu(\vartheta),\lambda(\vartheta))- N_{\alpha}(x(\vartheta))^TK_{\alpha}(x(\vartheta))^{-1}N_{\alpha}(x(\vartheta)).
 \end{array}
$$
To do this, define $A\in \Re^{m\times m}$ by
$$
A(\vartheta)=\left[
\begin{array}{l}
{\cal J}_yh(x(\vartheta),y(\vartheta))\\[3pt]
{\cal J}_yg_{\alpha}(x(\vartheta),y(\vartheta))\\[3pt]
\bar A
\end{array}
\right],
$$
where $\bar A\in \Re^{m-m_1-|\alpha|}$ is chosen such that $A(\vartheta_0)$ is nonsingular. Then if $\vartheta$ is close to $\vartheta_0$ enough, we have that $A(\vartheta)$ is nonsingular as well. By applying the standard Gram-Schmidt orthogonalization procedure to the columns of $A(\vartheta)$, we obtain an orthogonal matrix $P(\vartheta)=[P_1(\vartheta)\,\,P_2(\vartheta)]\in \Re^{m\times m}$ with $P_1(\vartheta)\in \Re^{m\times m_1+|\alpha|},P_2(\vartheta)\in \Re^{m-m_1-|\alpha|}$. Then $P(\vartheta)$ is a continuous function over $\textbf{B}(\vartheta_0,\delta)$. Let $Z(\vartheta)=P_2(\vartheta)$ satisfy (\ref{ta1}) and $Z(\vartheta)$ be  continuous  over $\textbf{B}(\vartheta_0,\delta)$. Then $Z(\vartheta)$ satisfies (\ref{ta1}) and (\ref{ta2}) when $\vartheta \in \textbf{B}(\vartheta_0,\delta)$ for small $\delta>0$, which comes from the fact that
$$
 Z(\vartheta_0)^T\Psi(\vartheta_0)Z(\vartheta_0) \succ 0,
$$
 from (\ref{eq:secSCs}).  The proof is completed. \hfill $\Box$

\section{Strong Regularity without Strict Complementarity}
\setcounter{equation}{0}
\quad \, In the Jacobian uniqueness condition of Problem (\ref{cminimax})  by Definition \ref{def:CS}, a critical condition is the strict complementarity for the upper level problem. In this section, we consider the case when this condition does not hold.

  Let $(x^*,u^*,v^*,y^*,\mu^*,\lambda^*)$ be a KKT point of Problem (\ref{cminimax}); namely, it
satisfies the following conditions
\begin{equation}\label{KKT}
\begin{array}{l}
\nabla_x {\cal L}(x,y,\mu,\lambda)+{\cal J}H(x)^Tu+{\cal J}G(x)^Tv=0,\\[3pt]
H(x)=0,\quad
0\leq v \perp G(x) \leq 0,\\[3pt]
\nabla_y {\cal L}(x,y,\mu,\lambda)=0,\\[3pt]
h(x,y)=0,\quad
0\geq \lambda \perp g(x,y) \leq 0.
\end{array}
\end{equation}
Let $z:=(x,u,v,y,\mu,\lambda)$ and define
$$
{\cal K}=\Re^n \times \Re^{n_1} \times \Re^{n_2}_+\times \times \Re^m\times \Re^{m_1}\times \Re^{m_2}_+
$$
and
\begin{equation}\label{GEf}
{\cal H}(z)=\left
[
\begin{array}{c}
 \nabla_x {\cal L}(x,y,\mu,\lambda)+{\cal J}H(x)^Tu+{\cal J}G(x)^Tv\\[3pt]
H(x)\\[3pt]
-G(x)\\[3pt]
\nabla_y {\cal L}(x,y,\mu,\lambda)\\[3pt]
h(x,y)\\[3pt]
-g(x,y)
\end{array}
\right].
\end{equation}
The KKT conditions above can be expressed as the following generalized equation
\begin{equation}\label{GeKKT}
0\in {\cal H}(z)
+N_{\cal K}(z).
\end{equation}
For $\eta=(\eta_x;\eta_H;\eta_G;\eta_y;\eta_h;\eta_g)$, it is easy to see that the perturbed generalized equation
\begin{equation}\label{pGeKKT}
\eta\in {\cal H}(z)
+N_{\cal K}(z)
\end{equation}
represents the KKT conditions for the  following canonical perturbation of
 Problem (\ref{cminimax}),
\begin{equation}\label{Cpcminimax}
\min_{x \in \Phi (\eta_H,\eta_G)}\max_{y \in Y(x,\eta_h,\eta_g)}f(x,y)-\langle \eta_x,x\rangle-\langle\eta_y,y \rangle,
\end{equation}
where $f:\Re^n\times \Re^m \rightarrow \Re$, $\Phi:\Re^{n_1}\times \Re^{n_2}\rightrightarrows \Re^n$ is a set-valued mapping defined by
\begin{equation}\label{CpPhi}
\Phi (\eta_H,\eta_G)=\Big\{x \in \Re^n: H(x)-\eta_H=0, G(x)+\eta_G\leq 0\Big\}
\end{equation}
and $Y: \Re^n \times \Re^{m_1}\times \Re^{m_2}\rightrightarrows \Re^m$ is a set-valued mapping defined by
\begin{equation}\label{CpYx}
Y(x,\eta_h,\eta_g)=\Big\{y \in \Re^m: h(x,y)-\eta_h=0, g(x,y)+\eta_g\leq 0\Big\}.
\end{equation}
Robinson \cite{Robinson80} introduced the  concept of strong regularity for a solution of the generalized
equation (\ref{GeKKT}).
\begin{definition}\label{Def.3.4}
 Let $z^*$ be a solution of the generalized equation (\ref{GeKKT}). We say that $z^*$  is a strongly regular
solution of the generalized equation (\ref{GeKKT}) if there exist positive numbers $\delta$ and $\varepsilon>0$ such that for every $\eta \in \textbf{B}(0, \delta)$, the following linearized generalized equation
\begin{equation}\label{LGeKKT}
\eta\in {\cal J}{\cal H}(z^*)(z-z^*)
+N_{\cal K}(z)
\end{equation}
has a unique solution in $\textbf{B}(z^*,\varepsilon)$, denoted by $\widehat z(\eta)$, and the mapping
$\widehat z: \textbf{B}(0, \delta)  \rightarrow \textbf{B}(z^*,\varepsilon)$ is Lipschitz continuous.
\end{definition}
It follows from \cite{Robinson80} or \cite{BS00} that if $z^*$  is a strongly regular
solution of the generalized equation (\ref{GeKKT}), then there exist positive numbers $\delta$ and $\varepsilon>0$
 such that for every $\eta \in \textbf{B}(0, \delta)$, the following generalized equation
\begin{equation}\label{LGeKKT}
\eta\in {\cal H}(z)
+N_{\cal K}(z)
\end{equation}
has a unique solution in $\textbf{B}(z^*,\varepsilon)$, denoted by $ z(\eta)$, and the mapping
$z: \textbf{B}(0, \delta)  \rightarrow \textbf{B}(z^*,\varepsilon)$ is Lipschitz continuous over $\textbf{B}(0, \delta)$.

To study the strong regularity of the KKT system at $(x^*,u^*,v^*,y^*,\mu^*,\lambda^*)$, we introduce the following definition.
\begin{definition}\label{def:wCS}
We say that Problem (\ref{cminimax}) satisfies  {\sc Property A} at $(x^*,u^*,v^*,y^*,\mu^*,\lambda^*)\in \Re^n\times \Re^{n_1}\times \Re^{n_2}\times \Re^m \times \Re^{m_1}\times \Re^{m_2}$ if
\begin{itemize}
\item[{\rm (i)}] $x^* \in \Phi$ and conditions in (\ref{KKT-P}) are satisfied at $(x^*,u^*,v^*,y^*,\mu^*,\lambda^*)$.
\item[{\rm (ii)}]The set vectors $\big\{\nabla H_1(x^*),\ldots,\nabla H_{n_1}\big\}\cup \Big\{\nabla G_i(x^*):i\in I(x^*)\Big\}$
are linearly independent, where $I(x^*)=\{i: G_i(x^*)=0, i=1.\ldots, n_2\}$.
\item[{\rm (iv)}] $y^* \in Y(x^*)$ and  Problem $({\rm P}_{x^*})$  satisfies  Jacobian uniqueness conditions at $(y^*,\mu^*,\lambda^*)$.
    \item[{\rm (v)}]For every $d_x \in {\rm Aff}\,{\cal C}(x^*)\setminus \{0\}$ (where ${\cal C}(x^*)$ is defined by (\ref{eq:CritiCS})),
         \begin{equation}\label{eq:StrsecSCs}
\begin{array}{l}
\displaystyle \sup_{(u,v) \in \Lambda(x^*)} \left\{ \left\langle \left[\displaystyle \sum_{j=1}^{n_1}u_i\nabla^2_{xx}H_j(x^*)+\displaystyle \sum_{i=1}^{n_2} v_i\nabla^2_{xx}G_i(x^*)\right]d_x,d_x \right\rangle\right\}\\[16pt]
 \quad \quad \quad \, +\left\langle \left[ \nabla^2_{xx}{\cal L}(x^*,y^*,\mu^*,\lambda^*)- N_{\alpha}(x^*)^TK_{\alpha}(x^*)^{-1}N_{\alpha}(x^*)\right]d_x,d_x \right\rangle > 0,
 \end{array}
\end{equation}
where $K_{\alpha}(x)$ is defined by (\ref{eq:Kx}) and $N_{\alpha}(x)$ is defined by (\ref{Nx}).
\end{itemize}
\end{definition}
It can be checked that the perturbed  Kojima mapping of the form
\begin{equation}\label{PKojima-f}
F(x,u,w,y,\mu,\xi)=\left[
\begin{array}{c}
\nabla_x{\cal L}(x,y,\mu,\xi^+)+{\cal J}H(x)^Tu+{\cal J}G(x)^Tw^+\\[3pt]
H(x)\\[3pt]
G(x)-w^-\\[3pt]
\nabla_y{\cal L}(x,y,\mu,\xi^+)\\[3pt]
h(x,y)\\[3pt]
-g(x,y)+\xi^-
\end{array}
\right
]-\left[
\begin{array}{c}
\eta_x\\[3pt]
\eta_H\\[3pt]
-\eta_G\\[3pt]
\eta_y\\[3pt]
\eta_h\\[3pt]
\eta_g
\end{array}
\right]
\end{equation}
is the Kojima mapping of the   canonical perturbation
 Problem (\ref{Cpcminimax}),
where $\Phi:\Re^{n_1}\times \Re^{n_2}\rightrightarrows \Re^n$ is  defined by
(\ref{CpPhi})
and $Y: \Re^n \times \Re^{m_1}\times \Re^{m_2}\rightrightarrows \Re^m$ is defined by
(\ref{CpYx}). Thus we have that $F(x,u,w,y,\mu,\xi)=\eta$ with $\eta=(\eta_x;\eta_H;-\eta_G;\eta_y;\eta_h;\eta_g)$ corresponds to KKT conditions for Problem (\ref{Cpcminimax}).
\begin{lemma}\label{lem:Snsing}
Suppose that {\sc Property A} of Problem (\ref{cminimax}) is satisfied at $(x^*,u^*,v^*,y^*,\mu^*,\lambda^*)\in \Re^n\times \Re^{n_1}\times \Re^{n_2}\times \Re^m \times \Re^{m_1}\times \Re^{m_2}$. Then for $w^*=v^*+G(x^*)$ and
$\xi^*=\lambda^*+g(x^*,y^*)$, any element of  $\partial F(x^*,u^*,w^*,y^*,\mu^*,\xi^*)$ is nonsingular.
\end{lemma}
{\bf Proof}. Since $\lambda^*-g(x^*,y^*)>0$, we know that $\xi^+$ and $\xi^-$ are differentiable at $\xi^*$,
and $w^+$ and $w^-$  are differentiable at $w^*$. Since  $w^+$ and $w^-$ are strongly semi-smooth in the sense of \cite{Qi93a},  $F$ is strongly semi-smooth at $(x^*,y^*,u^*,w^*,\mu^*,\xi^*)$. Without loss of generality, we assume
that
$$
\beta:=I(x^*)=\{1,\ldots, r\},\quad \alpha=\{1,\ldots, s\},
$$
where
$$
I(x^*)=\{i:G_i(x^*)=0, i=1,\ldots, n_2\},\quad \alpha=\{i: g_i(x^*,y^*)=0,i=1,\ldots,m_2\}.
$$
Let
$$
\beta_+=\Big\{i\in \beta: v^*_i>0\Big\},\quad \beta_0=\Big\{i\in \beta: v^*_i=0\Big\}
$$
and assume
$$
\beta_+=\Big\{1,\ldots, r_1\Big\}.
$$
Then, for $\beta^c=\{1,\ldots,n_2\}\setminus \beta$ and $\alpha^c=\{1,\ldots,m_2\}\setminus \alpha$, we get that
$$
\beta_+=\Big\{1,\ldots, r_1\Big\},\quad \beta_0=\Big\{r_1+1,\ldots, r\Big\},\quad \{\beta^c=\{r+1,\ldots, n_2\}, \quad \alpha^c=\{s+1,\ldots, m_2\}.
$$
Thus we obtain
\begin{equation}\label{SJxi}
{\cal J}\xi^+|_{\xi=\xi^*}=\left[
\begin{array}{cc}
I_s & 0\\[3pt]
0 & 0
\end{array}
\right],\quad
 {\cal J}\xi^-|_{\xi=\xi^*}=\left[
\begin{array}{cc}
0 & 0\\[3pt]
0 & I_{m_2-s}
\end{array}
\right]
\end{equation}
 and
\begin{equation}\label{eq:SgJ}
\begin{array}{l}
\partial w^+|_{w=w^*}=\left\{\left[
\begin{array}{ccc}
I_{r_1}&0 & 0\\[3pt]
0 & \omega_{\beta_0}& 0\\[3pt]
0 & 0& 0
\end{array}
\right]:\omega_{\beta_0}={\rm Diag}[\omega_{r_1+1},\cdot,\omega_r],\ \omega_i \in [0,1],\, i \in \beta_0\right\},\\[18pt]
\partial w^-|_{w=w^*}=\Big\{I_{n_2}-\omega: \omega\in \partial w^+|_{w=w^*}\Big\}.
\end{array}
\end{equation}

Denote
$$
\begin{array}{ll}
G^*_{11}=\nabla^2_{xx}{\cal L}(x^*,y^*,\mu^*,\lambda^*)
+\displaystyle \sum_{j=1}^{n_1}u_i\nabla^2_{xx}H_j(x^*)+\displaystyle \sum_{i=1}^{n_2} v_i\nabla^2_{xx}G_i(x^*),\\[4pt]
G^*_{12}=\nabla^2_{xy}{\cal L}(x^*,y^*,\mu^*,\lambda^*),\,\,
G^*_{22}=\nabla^2_{yy}{\cal L}(x^*,y^*,\mu^*,\lambda^*).
 \end{array}
$$
For simplicity, we use notations ${\cal J}_xh^*$ and ${\cal J}_yh^*$ to represent ${\cal J}_xh(x^*,y^*)$ and ${\cal J}_yh(x^*,y^*)$, respectively. The same
notations are also applied to $g_{\alpha}$ and $g_{\alpha^c}$. Let $V$ be an element of
$\partial F(x^*,u^*,w^*,y^*,\mu^*,\xi^*)$. Then there exists an element $\omega \in  \partial w^+|_{w=w^*}$
such that
{\small
\begin{equation}\label{SV}
\begin{array}{l}
V
=\left[
\begin{array}{ccccccccc}
G^*_{11}& {\cal J}H(x^*)^T & {\cal J}G_{\beta_+}(x^*)^T &\omega_{\beta_0} {\cal J}G_{\beta_0}(x^*)^T &0 &G^*_{12}&{\cal J}_xh^{*T}&-{\cal J}_xg_{\alpha}^{*T} &0\\[6pt]
{\cal J}H(x^*) & 0 &0 &0 &0 &0 &0 &0 &0 \\[6pt]
{\cal J}G_{\beta_+}(x^*) &0& 0&0&0&0&0&0&0\\[6pt]
{\cal J}G_{\beta_0}(x^*) &0& 0&-I_{|\beta_0|}+\omega_{\beta_0}&0&0&0&0&0\\[6pt]
{\cal J}G_{\beta^c}(x^*) & 0&0&0&-I_{n_2-r}&0&0&0&0\\[6pt]
G^{*T}_{12} & 0&0&0&0&G^*_{22}&{\cal J}_yh^{*T}&-{\cal J}_yg_{\alpha}^{*T} &0\\[6pt]
{\cal J}_xh^* & 0&0&0&0&{\cal J}_yh^*&0&0&0\\[6pt]
-{\cal J}_xg_{\alpha}^* & 0&0&0&0&-{\cal J}_yg_{\alpha}^*&0&0&0\\[6pt]
-{\cal J}_xg_{\alpha^c}^* &0& 0&0&0&-{\cal J}_yg_{\alpha^c}^*&0&0&I_{m_2-s}
\end{array}
\right].
\end{array}
\end{equation}}
The nonsingularity of $V$ is equivalent to the nonsingularity of the following matrix
{\small
$$
\begin{array}{l}
\left[
\begin{array}{ccccccccc}
G^*_{11}& {\cal J}H(x^*)^T & {\cal J}G_{\beta_+}(x^*)^T & \omega_{\beta_0}{\cal J}G_{\beta_0}(x^*)^T&G^*_{12}&{\cal J}_xh^{*T}&-{\cal J}_xg_{\alpha}^{*T} &0&0 \\[6pt]
{\cal J}H(x^*) & 0&0 &0 &0 &0 &0 &0 &0 \\[6pt]
{\cal J}G_{\beta_+}(x^*) &0& 0&0&0&0&0&0&0\\[6pt]
{\cal J}G_{\beta_0}(x^*) &0& 0&-I_{|\beta_0|}+\omega_{\beta_0}&0&0&0&0&0\\[6pt]
G^{*T}_{12} & 0&0&0&G^*_{22}&{\cal J}_yh^{*T}&-{\cal J}_yg_{\alpha}^{*T} &0&0\\[6pt]
{\cal J}_xh^* &0& 0&0&{\cal J}_yh^*&0&0&0&0\\[6pt]
-{\cal J}_xg_{\alpha}^*&0 & 0&0&-{\cal J}_yg_{\alpha}^*&0&0&0&0\\[6pt]
{\cal J}G_{\beta^c}(x^*)&0 & 0&0&0&0&0&-I_{n_2-r}&0\\[6pt]
-{\cal J}_xg_{\alpha^c}^*&0 & 0&0&-{\cal J}_yg_{\alpha^c}^*&0&0&0&I_{m_2-s}
\end{array}
\right],
\end{array}
$$}
which is equivalent to the nonsingularity of the following
matrix
\begin{equation}\label{SJFa}
H(\omega)=
\left[
\begin{array}{ccccccc}
G^*_{11}& {\cal J}H(x^*)^T &  {\cal J}G_{\beta_+}(x^*)^T & \omega_{\beta_0}{\cal J}G_{\beta_0}(x^*)^T &G^*_{12}&{\cal J}_xh^{*T}&-{\cal J}_xg_{\alpha}^{*T}  \\[6pt]
{\cal J}H(x^*) &0 & 0 &0 &0&0&0  \\[6pt]
{\cal J}G_{\beta_+}(x^*)&0 & 0&0&0&0&0\\[6pt]
{\cal J}G_{\beta_0}(x^*)&0 & 0&-I_{|\beta_0|}+\omega_{\beta_0}&0&0&0\\[6pt]
G^{*T}_{12} &0& 0&0&G^*_{22}&{\cal J}_yh^{*T}&-{\cal J}_yg_{\alpha}^{*T}\\[6pt]
{\cal J}_xh^*&0 & 0&0&{\cal J}_yh^*&0&0\\[6pt]
-{\cal J}_xg_{\alpha}^*&0 & 0&0&-{\cal J}_yg_{\alpha}^*&0&0
\end{array}
\right].
\end{equation}
Therefore, we only need to prove that the matrix $H(\omega)$ is nonsingular. From Lemma \ref{ns}, we obtain
that $K_{\alpha}(x^*)$ is nonsingular, where
$$
K_{\alpha}(x^*)=\left[
\begin{array}{ccc}
G^*_{22}&{\cal J}_yh^{*T}&-{\cal J}_yg_{\alpha}^{*T}\\[6pt]
{\cal J}_yh^*&0&0\\[6pt]
-{\cal J}_yg_{\alpha}^*&0&0
\end{array}
\right].
$$
So we only need to prove that
$
H(\omega)/K_{\alpha}(x^*)
$ is nonsingular. Notice that
\begin{equation}\label{SHkf}
\begin{array}{rcl}
H(\omega)/K_{\alpha}(x^*)&=& \left[
\begin{array}{cccc}
G^*_{11}& {\cal J}H(x^*)^T & {\cal J}G_{\beta_+}(x^*)^T & \omega_{\beta_0}{\cal J}G_{\beta_0}(x^*)^T \\[6pt]
{\cal J}H(x^*) & 0 &0&0 \\[6pt]
{\cal J}G_{\beta_+}(x^*) & 0&0&0\\[6pt]
{\cal J}G_{\beta_0}(x^*) & 0&0&-I_{|\beta_0|}+\omega_{\beta_0}\\[6pt]
\end{array}
\right]\\[20pt]
&&-\left[
\begin{array}{ccc}
G^*_{12}&{\cal J}_xh^{*T}&-{\cal J}_xg_{\alpha}^{*T}  \\[3pt]
0 &0&0  \\[3pt]
0&0&0\\[3pt]
0&0&0
\end{array}
\right]K_{\alpha}(x^*)^{-1}\left[
\begin{array}{cccc}
G^{*T}_{12}&0&0 &0 \\[6pt]
{\cal J}_xh^* &0&0&0  \\[6pt]
-{\cal J}_xg_{\alpha}^*&0&0&0
\end{array}
\right]\\[20pt]
&=& \left[
\begin{array}{cccc}
G^*_{11}-N_{\alpha}(x^*)^TK_{\alpha}(x^*)N_{\alpha}(x^*)& {\cal J}H(x^*)^T & {\cal J}G_{\beta_+}(x^*)^T & \omega_{\beta_0}{\cal J}G_{\beta_0}(x^*)^T \\[6pt]
{\cal J}H(x^*) & 0 &0&0 \\[6pt]
{\cal J}G_{\beta_+}(x^*) & 0&0&0\\[6pt]
{\cal J}G_{\beta_0}(x^*) & 0&0&-I_{|\beta_0|}+\omega_{\beta_0}\\[6pt]
\end{array}
\right].
\end{array}
\end{equation}
It is easy to check  that ${\rm Aff}\,{\cal C}(x^*)$ is of the following subspace
\begin{equation}\label{eq:Skerhg}
{\rm Aff}\,{\cal C}(x^*)={\rm Ker}\, \left[
\begin{array}{l}
{\cal J}H(x^*) \\[6pt]
{\cal J}G_{\beta_+}(x^*)
\end{array}
\right
].
\end{equation}
Now we prove that $H(\omega)/K_{\alpha}(x^*)$ is nonsingular via the formula (\ref{SHkf}). Let $a_1\in \Re^n$, $a_2\in \Re^{n_1}$, $a_3\in \Re^{r_1}$ and $a_4\in \Re^{r-r_1}$ satisfy
$$
\left[
\begin{array}{cccc}
G^*_{11}-N_{\alpha}(x^*)^TK_{\alpha}(x^*)N_{\alpha}(x^*)& {\cal J}H(x^*)^T & {\cal J}G_{\beta_+}(x^*)^T & \omega_{\beta_0}{\cal J}G_{\beta_0}(x^*)^T \\[6pt]
{\cal J}H(x^*) & 0 &0&0 \\[6pt]
{\cal J}G_{\beta_+}(x^*) & 0&0&0\\[6pt]
{\cal J}G_{\beta_0}(x^*) & 0&0&-I_{|\beta_0|}+\omega_{\beta_0}\\[6pt]
\end{array}
\right]\left[
\begin{array}{c}
a_1\\[6pt]
a_2\\[6pt]
a_3\\[6pt]
a_4
\end{array}
\right
]=0_{n+n_1+r}
$$
or
\begin{eqnarray}
&&\hspace{-1.3cm}[G^*_{11}-N_{\alpha}(x^*)^TK_{\alpha}(x^*)N_{\alpha}(x^*)]a_1+{\cal J}H(x^*)^Ta_2
+{\cal J}G_{\beta_+}(x^*)^Ta_3
+ \omega_{\beta_0}{\cal J}G_{\beta_0}(x^*)^Ta_4 =0,
\label{Sabc}\\[6pt]
&&\hspace{-1.3cm}{\cal J}H(x^*)a_1=0,\label{Sabc-2}\\[6pt]
&&\hspace{-1.3cm}{\cal J}G_{\beta_+}(x^*)a_1=0,\label{Sabc-3}\\[6pt]
&&\hspace{-1.3cm}{\cal J}G_{\beta_0}(x^*)a_1
+ [-I_{|\beta_0|}+\omega_{\beta_0}]a_4=0.\label{Sabc-4}
\end{eqnarray}
It follows from (\ref{Sabc-2}) and (\ref{Sabc-3}) that $a_1 \in {\rm Aff}\,{\cal C}(x^*)$.
Premultiplying $a^T_1$ to (\ref{Sabc}), we obtain
\begin{equation}\label{a1times}
a_1^T[G^*_{11}-N_{\alpha}(x^*)^TK_{\alpha}(x^*)N_{\alpha}(x^*)]a_1+a_1^T \omega_{\beta_0}{\cal J}G_{\beta_0}(x^*)^Ta_4=0.
\end{equation}
From the relation (\ref{Sabc-4}), we get that
\begin{equation}\label{4e}
a_1^T \omega_{\beta_0}{\cal J}G_{\beta_0}(x^*)^Ta_4
=\displaystyle \sum_{i\in \beta_0:0< w_i<1}\displaystyle \frac{\omega_i}{1-\omega_i}[\nabla G_i(x^*)^Ta_1]^2\geq0.
\end{equation}
If follows from  (\ref{Sabc}) and (\ref{4e}) that
$$
\begin{array}{l}
\displaystyle  \left\langle \left[\displaystyle \sum_{j=1}^{n_1}u^*_i\nabla^2_{xx}H_j(x^*)+\displaystyle \sum_{i=1}^{n_2} v^*_i\nabla^2_{xx}G_i(x^*)\right]a_1,a_1 \right\rangle\\[16pt]
 \quad \quad \quad \, +\left\langle \left[ \nabla^2_{xx}{\cal L}(x^*,y^*,\mu^*,\lambda^*)- N_{\alpha}(x^*)^TK_{\alpha}(x^*)^{-1}N_{\alpha}(x^*)\right]a_1,a_1 \right\rangle\leq 0,
 \end{array}
$$
which implies $a_1=0$ from the condition (\ref{eq:secSCs}). Let $\bar \beta_0=\{i\in \beta_0:\omega_i=1\}$ and $\bar \beta_0^c=
\beta_0\setminus \bar \beta_0$. Then from
(\ref{Sabc-4}), $a_1=0$ and (\ref{Sabc}), we can get that
\begin{equation}\label{a4c}
[a_4]_{\bar \beta_0^c}=0
\end{equation}
and
\begin{equation}\label{a4cb}
{\cal J}H(x^*)^T a_2+ {\cal J}G_{\beta_+}(x^*)^Ta_3+{\cal J}G_{\bar \beta_0}(x^*)^T[a_4]_{\bar \beta_0}=0.
\end{equation}
From
 (iii), we obtain $a_2=0$,$a_3=0$ and $[a_4]_{\bar \beta_0}=0$. Combining with  (\ref{a4c}), we have that $a_2=0$, $a_3=0$ and $a_4=0$. Therefore $H(\omega)/K_{\alpha}(x^*)$ is nonsingular.  The proof is completed. \hfill $\Box$

 Now we are in a position to establish the main result about the strong regularity of the KKT system for
Problem (\ref{cminimax}).
\begin{theorem}\label{th:sR}
Let $(x^*,y^*)\in \Re^n \times \Re^m$ be a point around which $f$, $h$ and $g$ are twice differentiable, and $H$ and $G$ are twice continuously differentiable. Assume that there exists  $(u^*,v^*,\mu^*,\lambda^*)\in  \Re^{n_1}\times \Re^{n_2}\times \Re^{m_1}\times \Re^{m_2}$ such that $(x^*,u^*,v^*,y^*,\mu^*,\lambda^*)$ satisfies Karush-Kuhn-Tucker conditions for
Problem (\ref{cminimax}).  Consider
the following four statements
\begin{itemize}
\item[{\rm (a)}]  {\sc Property A} holds at $(x^*,u^*,v^*,y^*,\mu^*,\lambda^*)$.
\item[{\rm (b)}] For $w^*=v^*+G(x^*)$ and
$\xi^*=\lambda^*+g(x^*,y^*)$, any element of  $\partial F(x^*,u^*,w^*,y^*,\mu^*,\xi^*)$ is nonsingular.
\item[{\rm (c)}] $F$ is  a locally Lipschitz homeomorphism near $(x^*,u^*,w^*,y^*,\mu^*,\xi^*)$.
\item[{\rm (d)}] The  point $(x^*,u^*,v^*,y^*,\mu^*,\lambda^*)$  is a strongly regular solution of the generalized equation (\ref{GeKKT}).
\end{itemize}
Then it holds that (a) $\Longrightarrow$ (b) $\Longrightarrow$ (c) $\Longleftrightarrow$ (d).
\end{theorem}
{\bf Proof}. From Lemma \ref{lem:Snsing}, we obtain (a)$\Longrightarrow$ (b).  By Clarke's inverse function theorem (Clarke \cite{Clarke76,Clarke83}), $F$ is a locally Lipschitz homeomorphism near $(x^*,u^*,w^*,y^*,\mu^*,\xi^*)$ and hence we get (b)$\Longrightarrow$ (c).
Noting that the generalized equation
$$
\eta\in {\cal H}(z)
+N_{\cal K}(z)
$$
represents the KKT conditions for
Problem (\ref{Cpcminimax}) and $F(x,u,w,y,\mu,\xi)=\eta$ with $\eta=(\eta_x$; $\eta_H$; $-\eta_G$; $\eta_y$; $\eta_h$; $\eta_g)$ corresponds to KKT conditions for Problem (\ref{Cpcminimax}), we obtain the equivalence between (c) and (d). The proof is completed.
\hfill $\Box$
 \section{Some Concluding Remarks}\label{Sec4}
In this paper, we have analyzed the stability properties of the Kurash-Kuhn-Tucker (KKT) system for Problem (\ref{cminimax}). Firstly, we proposed the definition of Jacobian uniqueness condition of Problem (\ref{cminimax}) and proved that this property is stable with respect to a small ${\cal C}^2$-perturbation.  Secondly, comparing with the Jacobian uniqueness condition, we proposed {\sc Property A} by eliminating the strict complementarity condition for the outer level constraints and adopting the strong second-order sufficiency optimality condition. We proved that the strong regularity of  the KKT system at the KKT point is equivalent to the  local Lipschitz homeomorphism of the Kojima mapping near the KKT point. Finally, we proved that
{\sc Property A} is a sufficient condition for the strong regularity of  the KKT system at the KKT point.

 There are many problems about the stability of constrained minimax optimization left to us. For instance, in our analysis, even {\sc Property A} requires the Jacobian uniqueness condition for the inner level problem. Is it possible to weaken this condition? A closely related problem is how to obtain the second-order optimality conditions for the constrained minimax problem when
 the Jacobian uniqueness condition for the inner level problem fails. This might be a difficult problem. Another question is, when the constraints are linear and the objective function is even convex-concave, can we have a sharp theoretical result like the result for linear semidefinte programming in \cite{CS08}?

%=========References=============\\
 
\end{document}